\documentclass[11pt]{article}
\usepackage{amssymb}

\usepackage{amsmath}
\usepackage{theorem}
\usepackage{here}
\usepackage[dvipdfmx]{color}

\newtheorem{thm}{Theorem}[section]
\newtheorem{lem}{Lemma}[section]

\theorembodyfont{\rmfamily}

\newtheorem{remark}{Remark}[section]
\makeatletter

\@addtoreset{equation}{section}
\makeatother

\usepackage{comment} % 
\usepackage{bm}
\usepackage[dvipdfmx]{graphicx} % 

\def\pt{{\tilde p}}
\def\nt{{\tilde n}}

\def\pd{\partial}

\def\bvth{{\text{\boldmath $\vartheta$}}}
\def\Th{{\Theta}}
\def\De{{\Delta}}

\def\Ga{{\Gamma}}

\def\d{{\text{\boldmath $d$}}}

\def\j{{\text{\boldmath $j$}}}

\def\p{{\text{\boldmath $p$}}}

\def\x{{\text{\boldmath $x$}}}

\def\X{{\text{\boldmath $X$}}}
\def\Y{{\text{\boldmath $Y$}}}

\def\ph{{\hat p}}

\def\pbh{{\hat \p}}

\def\pbt{{\tilde \p}}

\def\diag{{\rm diag\,}}

\def\[{{\text{\boldmath $[$}}}
\def\]{{\text{\boldmath $]$}}}

\def\/{{\Bigr/\!\!}}

\def\1r{{\rm (1)}}
\def\2r{{\rm (2)}}
\def\3r{{\rm (3)}}
\def\4r{{\rm (4)}}
\def\5r{{\rm (5)}}

\def\non{{\nonumber}}

\begin{document}
\title{Bayesian Estimation of Multinomial Cell Probabilities Incorporating Information from Aggregated Observations\footnote{This preprint has not undergone peer review (when applicable) or any post-submission improvements or corrections. 
The Version of Record of this article is published in Japanese Journal of Statistics and Data Science, and is available online at https://doi.org/10.1007/s42081-023-00224-z. }}
\author{
Yasuyuki Hamura\footnote{Graduate School of Economics, Kyoto University, 
Yoshida-Honmachi, Sakyo-ku, Kyoto, 606-8501, JAPAN. 
\newline{
E-Mail: yasu.stat@gmail.com}} \
}
\maketitle
\begin{abstract}
In this note, we consider the problem of estimating multinomial cell probabilities under the entropy loss when side information in aggregated data is available. 
We use the Jeffreys prior to obtain Bayes estimators. 
It is shown that by incorporating the side information, we can construct an improved estimator. 

\par\vspace{4mm}
{\it Key words and phrases:\ aggregated data, dominance, entropy loss, multinomial distribution. } 
\end{abstract}

\section{Introduction}
\label{sec:introduction}
Simultaneous estimation of multinomial cell probabilities has been studied by many authors. 
For example, Fienberg and Holland (1973) proposed smoothed estimators based on geometrical and Bayesian arguments. 
Leonard (1977), Albert (1981), Albert and Gupta (1983), Alam and Mitra (1986), and Albert (1987) considered different Bayesian estimators. 

Albert (1985) and Gibbons and Greenberg (1994) considered problems where there are partially classified data; they considered contingency tables when some data are classified only by rows (columns) and not by columns (rows). 
However, it seems that theoretical support has not been fully provided, as far as the author knows. 
In this paper, we treat a similar problem in a decision-theoretic framework. 

Specifically, we consider the case where we observe (1) fully classified data and (2) additional data which are classified by rows. 
We use the entropy loss for the multinomial distribution, under which it is unclear whether a minimax estimator can be derived. 
(For the quadratic loss, Nayak and Naik (1989) obtained a minimax Bayes estimator.) 
However, the entropy loss is directly related to predictive density estimation under the Kullback-Leibler divergence (see, for example, Komaki (2012)), and we show that incorporating the additional information can be useful under the entropy loss. 
Although such a result was obtained by Hamura (2021b) for the Poisson case, in this paper we take a different approach to establishing dominance.

\section{Main Results}
\label{sec:results} 
\subsection{The problem}
\label{subsec:problem} 
Let $\p = (( p_{i, j} )_{j = 0}^{n} )_{i = 0}^{m}$ be $(1 + m) (1 + n)$ unknown positive probabilities which sum to $1$. 
Suppose that we observe the direct observations 
\begin{align}
\X = (( X_{i, j} )_{j = 0}^{n} )_{i = 0}^{m} \sim {\rm{Multin}}_{(1 + m) (1 + n) - 1} ( \x | N, \p ) = {N ! \over \prod_{i = 0}^{m} \prod_{j = 0}^{n} x_{i, j} !} \prod_{i = 0}^{m} \prod_{j = 0}^{n} {p_{i, j}}^{x_{i, j}} \text{,} \non 
\end{align}
$\x = (( x_{i, j} )_{j = 0}^{n} )_{i = 0}^{m} \in %
\big\{ (( \mathring{x} _{i, j} )_{j = 0}^{n} )_{i = 0}^{m} \in {\mathbb{N} _0}^{(1 + m) (1 + n)} \big| \sum_{i = 0}^{m} \sum_{j = 0}^{n} \mathring{x} _{i, j} = N \big\} $, where $\mathbb{N} _0 = \{ 0 \} \cup \mathbb{N} = \{ 0, 1, 2, \dotsc \} $. 
Suppose that $\Y ^{*} = (( Y_{i, j}^{*} )_{j = 0}^{n} )_{i = 0}^{m} \sim {\rm{Multin}}_{(1 + m) (1 + n) - 1} ( N' , \p )$ independently of $\X $ and that we also observe 
\begin{align}
\Y = ( Y_i )_{i = 0}^{m} = \Big( \sum_{j = 0}^{n} Y_{i, j}^{*} \Big) _{i = 0}^{m} \sim {\rm{Multin}}_{m} ( N' , \p _{\cdot } ) \text{,} \non 
\end{align}
where $\p _{\cdot } = ( p_{i, \cdot } )_{i = 0}^{m} = \big( \sum_{j = 0}^{n} p_{i, j} \big) _{i = 0}^{m}$. 
We consider the problem of estimating $\p $ under the entropy loss 
\begin{align}
L( \d , \p ) &= \sum_{i = 0}^{m} \sum_{j = 0}^{n} p_{i, j} \log {p_{i, j} \over d_{i, j}} \text{,} \non 
\end{align}
$\d = (( d_{i, j} )_{j = 0}^{n} )_{i = 0}^{m} \in \Th = \big\{ (( \mathring{p} _{i, j} )_{j = 0}^{n} )_{i = 0}^{m} \in (0, \infty )^{(1 + m) (1 + n)} \big| \sum_{i = 0}^{m} \sum_{j = 0}^{n} \mathring{p} _{i, j} = 1 \big\} $.

\subsection{Bayes estimators}
\label{subsec:estimators} 
Here, we use the Jeffreys prior to obtain Bayes estimators. 
The estimators based on $\X $ and $( \X , \Y )$ are given in Lemma \ref{lem:estimators}. 
Let $X_{i, \cdot } = \sum_{j = 0}^{n} X_{i, j}$ for $i = 0, 1, \dots m$. 
\begin{lem}
\label{lem:estimators} 
\hfill
\begin{itemize}
\item[{\rm{(i)}}]
The Bayes estimator with respect to the Jeffreys prior and the direct observations $\X $ are 
\begin{align}
\pbt = (( \pt _{i, j} )_{j = 0}^{n} )_{i = 0}^{m} = \Big( \Big( {X_{i, j} + 1 / 2 \over N + (1 + m) (1 + n) / 2} \Big) _{j = 0}^{n} \Big) _{i = 0}^{m} \text{.} \non 
\end{align}
\item[{\rm{(ii)}}]
The Bayes estimator with respect to the Jeffreys prior and the direct and indirect observations $( \X , \Y )$ are 
\begin{align}
\pbh = (( \ph _{i, j} )_{j = 0}^{n} )_{i = 0}^{m} = \Big( \Big( {X_{i, \cdot } + Y_i + (1 + n) / 2 \over N + N' + (1 + m) (1 + n) / 2} {X_{i, j} + 1 / 2 \over X_{i, \cdot } + (1 + n) / 2} \Big) _{j = 0}^{n} \Big) _{i = 0}^{m} \text{.} \non 
\end{align}
\end{itemize}
\end{lem}

In the multinomial case, the Jeffreys prior is proper, and hence the associated estimator is admissible. 
However, we note that in the present setting where we can use $\Y $ as well as $\X $, it does not follow that $\pbt $ is admissible; it is admissible in the class of estimators which are functions of $\X $ only.

\subsection{Risk comparisons}
\label{subsec:dominance} 
In this section, we compare the risk functions of $\pbt $ and $\pbh $. 
Let $\De _{\p } (N, N' ) = E_{\p } [ L( \pbh , \p ) ] - E_{\p } [ L( \pbt , \p ) ]$. 
For notational simplicity, let $p_i = p_{i, \cdot }$ and $X_i = X_{i, \cdot }$ for $i = 0, 1, \dots , m$. 
\begin{lem}
\label{lem:decomposition} 
The difference in risks between $\pbh $ and $\pbt $ is decomposed as follows: 
\begin{align}
\De _{\p } (N, N' ) &= E_{\p } \Big[ \log {N + N' + (1 + m) (1 + n) / 2 \over N + (1 + m) (1 + n) / 2} + \sum_{i = 0}^{m} p_i \log {X_i + (1 + n) / 2 \over X_i + Y_i + (1 + n) / 2} \Big] \non \\
&= \sum_{\nu = 1}^{N'} \De _{\p } (N + \nu - 1, 1) \text{.} \non 
\end{align}
\end{lem}

In Lemmas \ref{lem:RD} and \ref{lem:max_RD} below, we consider $\De _{\p } (N, 1)$ without loss of generality. 
Define the function $H \colon (0, 1) \to (0, \infty )$ by 
\begin{align}
H(p) &= %
p^2 E \Big[ - \log \Big\{ 1 - {1 \over X(p) + 1 + (1 + n) / 2} \Big\} \Big] \text{,} \quad p \in (0, 1) \text{,} \non 
\end{align}
where $X(p) \sim {\rm{Bin}} (N, p)$ for $p \in (0, 1)$. 
\begin{lem}
\label{lem:RD} 
We have 
\begin{align}
\De _{\p } (N, 1) &= \log \Big\{ 1 + {1 \over N + (1 + m) (1 + n) / 2} \Big\} - \sum_{i = 0}^{m} H( p_i ) \text{.} \non 
\end{align}
\end{lem}

We can prove the convexity of the function $H$. 
\begin{lem}
\label{lem:convexity} 
Suppose that $n \ge 3$. 
Then the function $H$ is convex. 
\end{lem}

Let $\p ^{*} = (( p_{i, j}^{*} )_{j = 0}^{n} )_{i = 0}^{m} \in \Th $ be such that $p_{i, j}^{*} = 1 / \{ (1 + m) (1 + n) \} $ for all $i = 0, 1, \dots , m$ and all $j = 0, 1, \dots , n$. 
Then $\sum_{j = 0}^{n} p_{i, j}^{*} = 1 / (1 + m)$ for all $i = 0, 1, \dots , m$. 
\begin{lem}
\label{lem:max_RD} 
Suppose that $n \ge 3$. 
Then the risk difference $\De _{\p } (N, 1)$ is maximized at $\p = \p ^{*}$: 
\begin{align}
\max_{\p \in \Th } \De _{\p } (N, 1) &= \De _{\p ^{*}} (N, 1) \text{.} \non 
\end{align}
\end{lem}

Now, we have the following results. 
\begin{thm}
\label{thm:N_Y} 
Fix $m \in \mathbb{N}$, $n \in \mathbb{N} \cap [3, \infty )$, and $N \in \mathbb{N}$. 
Assume that $N'$ is sufficiently large. 
Then $\pbt $ is dominated by $\pbh $. 
\end{thm}
\begin{thm}
\label{thm:N_X} 
Fix $m \in \mathbb{N}$, $n \in \mathbb{N} \cap [3, \infty )$, and $N' \in \mathbb{N}$. 
Assume that $N$ is sufficiently large. 
Then $\pbt $ is dominated by $\pbh $. 
\end{thm}

\begin{remark}
The maximum risk difference 
\begin{align}
\De _{\p ^{*}} (N, 1) &= \log \Big\{ 1 + {1 \over N + (1 + m) (1 + n) / 2} \Big\} \non \\
&\quad - {1 \over m + 1} E \Big[ - \log \Big\{ 1 - {1 \over X(1 / (1 + m)) + 1 + (1 + n) / 2} \Big\} \Big] \non 
\end{align}
is not necessarily easy to analytically evaluate. 
It can be shown that $\{ \pd / ( \pd n) \} \De _{\p ^{*}} (1, 1) < 0$, which implies that $\pbh $ does not dominate $\pbt $ when $N = 1$. 
It can be shown that $\lim_{m \to \infty } (1 + m) \De _{\p ^{*}} (N, 1) > 0$, which implies that $\pbh $ does not dominate $\pbt $ as $m \to \infty $. 
Finally, it can be shown that $\lim_{n \to \infty } n^3 \{ \pd / ( \pd n) \} \De _{\p ^{*}} (N, 1) < 0$, which implies that $\pbh $ does not dominate $\pbt $ as $n \to \infty $. 
\end{remark}

\section{Proofs}
\label{sec:proofs} 
Here, we prove the results of the previous section. 

\bigskip

\noindent
{\bf Proof of Lemma \ref{lem:estimators}.} \ \ Let $\rho _{i, j} = p_{i, j} / p_i$ for $j = 0, 1, \dots , n$ for $i = 0, 1, \dots , m$. 
The probability mass function of $( \X , \Y )$ is 
\begin{align}
f( \X , \Y | \bvth ) &= \Big\{ {N ! \over \prod_{i = 0}^{m} \prod_{j = 0}^{n} X_{i, j} !} \prod_{i = 0}^{m} \prod_{j = 0}^{n} ( p_i \rho _{i, j} )^{X_{i, j}} \Big\} {N' ! \over \prod_{i = 0}^{m} Y_i !} \prod_{i = 0}^{m} {p_i}^{Y_i} \non \\
&= {N ! \over \prod_{i = 0}^{m} \prod_{j = 0}^{n} X_{i, j} !} {N' ! \over \prod_{i = 0}^{m} Y_i !} \Big( \prod_{i = 0}^{m} {p_i}^{X_i + Y_i} \Big) \prod_{i = 0}^{m} \prod_{j = 0}^{n} {\rho _{i, j}}^{X_{i, j}} \text{,} \non 
\end{align}
where $\bvth = (( p_i )_{i = 1}^{m} , (( \rho _{i, j} )_{j = 1}^{n} )_{i = 0}^{m} )$. 
Then 
\begin{align}
&{\pd \over \pd p_i} \log f( \X , \Y | \bvth ) = {X_i + Y_i \over p_i} - {X_0 + Y_0 \over p_0} \non 
\end{align}
for $i = 1, \dots , m$ and 
\begin{align}
&{\pd \over \pd \rho _{i, j}} \log f( \X , \Y | \bvth ) = {X_{i, j} \over \rho _{i, j}} - {X_{i, 0} \over \rho _{i, 0}} \text{,} \non 
\end{align}
for $j = 1, \dots , n$ for $i = 0, 1, \dots , m$. 
Therefore, 
\begin{align}
&{\pd ^2 \over ( \pd p_i ) ( \pd p_{i'} )} \log f( \X , \Y | \bvth ) = \begin{cases} \displaystyle - {X_i + Y_i \over {p_i}^2} - {X_0 + Y_0 \over {p_0}^2} \text{,} & \text{if $i = i'$} \text{,} \\ \displaystyle - {X_0 + Y_0 \over {p_0}^2} \text{,} & \text{if $i \neq i'$} \text{,} \end{cases} \non 
\end{align}
for $i, i' = 1, \dots , m$ and 
\begin{align}
&{\pd ^2\over ( \pd \rho _{i, j} ) ( \pd \rho _{i, j'} )} \log f( \X , \Y | \bvth ) = \begin{cases} \displaystyle - {X_{i, j} \over {\rho _{i, j}}^2} - {X_{i, 0} \over {\rho _{i, 0}}^2} \text{,} & \text{if $j = j'$} \text{,} \\ \displaystyle - {X_{i, 0} \over {\rho _{i, 0}}^2} \text{,} & \text{if $j \neq j'$} \text{,} \end{cases} \non 
\end{align}
for $j, j' = 1, \dots , n$ for $i = 0, 1, \dots , m$. 
Therefore, the Jeffreys prior is 
\begin{align}
\pi ^{\rm{J}} ( \bvth ) d\bvth &\propto \sqrt{\Big| E_{\p } \Big[ \diag \Big( {X_1 + Y_1 \over {p_1}^2}, \dots , {X_m + Y_m \over {p_m}^2} \Big) + \j ^{(m)} ( \j ^{(m)} )^{\top } {X_0 + Y_0 \over {p_0}^2} \Big] \Big| } \non \\
&\quad \times \prod_{i = 0}^{m} \sqrt{\Big| E_{\p } \Big[ \diag \Big( {X_{i, 1} \over {\rho _{i, 1}}^2}, \dots , {X_{i, n} \over {\rho _{i, n}}^2} \Big) + \j ^{(n)} ( \j ^{(n)} )^{\top } {X_{i, 0} \over {\rho _{i, 0}}^2} \Big] \Big| } \non \\
&\propto \sqrt{\Big| \diag \Big( {1 \over p_1}, \dots , {1 \over p_m} \Big) + \j ^{(m)} ( \j ^{(m)} )^{\top } {1 \over p_0} \Big| } \non \\
&\quad \times \Big( \prod_{i = 0}^{m} {p_i}^{n / 2} \Big) \prod_{i = 0}^{m} \sqrt{\Big| \diag \Big( {1 \over \rho _{i, 1}}, \dots , {1 \over \rho _{i, n}} \Big) + \j ^{(n)} ( \j ^{(n)} )^{\top } {1 \over \rho _{i, 0}} \Big| } \non \\
&= \sqrt{\Big| \diag \Big( {1 \over p_1}, \dots , {1 \over p_m} \Big) \Big| \Big[ 1 + {( \j ^{(m)} )^{\top } \over \sqrt{p_0}} \Big\{ \diag \Big( {1 \over p_1}, \dots , {1 \over p_m} \Big) \Big\} ^{- 1} {\j ^{(m)} \over \sqrt{p_0}} \Big] } \non \\
&\quad \times \Big( \prod_{i = 0}^{m} {p_i}^{n / 2} \Big) \prod_{i = 0}^{m} \sqrt{\Big| \diag \Big( {1 \over \rho _{i, 1}}, \dots , {1 \over \rho _{i, n}} \Big) \Big| \Big[ 1 + {( \j ^{(n)} )^{\top } \over \sqrt{\rho _{i, 0}}} \Big\{ \diag \Big( {1 \over \rho _{i, 1}}, \dots , {1 \over \rho _{i, n}} \Big) \Big\} ^{- 1} {\j ^{(n)} \over \sqrt{\rho _{i, 0}}} \Big] } \non \\
&= \Big\{ \prod_{i = 0}^{m} {p_i}^{(1 + n) / 2 - 1} \Big\} \prod_{i = 0}^{m} \prod_{j = 0}^{n} {\rho _{i, j}}^{1 / 2 - 1} \text{.} \non 
\end{align}
Since the Bayes estimators under the entropy loss are the posterior means of $\p = (( p_i \rho _{i, j} )_{j = 0}^{n} )_{i = 0}^{m}$, the result follows. 
\hfill$\Box$

\bigskip

\noindent
{\bf Proof of Lemma \ref{lem:decomposition}.} \ \ By Lemma \ref{lem:estimators}, we have 
\begin{align}
\De _{\p } (N, N' ) %
&= E_{\p } \Big[ \sum_{i = 0}^{m} \sum_{j = 0}^{n} p_{i, j} \log {\pt _{i, j} \over \ph _{i, j}} \Big] \non \\
&= E_{\p } \Big[ \sum_{i = 0}^{m} \sum_{j = 0}^{n} p_{i, j} \log \Big\{ {N + N' + (1 + m) (1 + n) / 2 \over N + (1 + m) (1 + n) / 2} {X_i + (1 + n)  / 2 \over X_i + Y_i + (1 + n)  / 2} \Big\} \Big] \non \\
&= E_{\p } \Big[ \log {N + N' + (1 + m) (1 + n) / 2 \over N + (1 + m) (1 + n)  / 2} + \sum_{i = 0}^{m} p_i \log {X_i + (n + 1) / 2 \over X_i + Y_i + (1 + n) / 2} \Big] \non \\
&= \sum_{\nu = 1}^{N'} \De _{\p } (N + \nu - 1, 1) \text{,} \non 
\end{align}
which is the desired result. 
\hfill$\Box$

\bigskip

\noindent
{\bf Proof of Lemma \ref{lem:RD}.} \ \ Let $( Z_i )_{i = 0}^{m} \sim {\rm{Multin}}_{m} (1, ( p_i )_{i = 0}^{m} )$ be an independent set of multinomial variables. 
Then, by Lemma \ref{lem:decomposition}, 
\begin{align}
\De _{\p } (N, 1) &= \log {N + 1 + (1 + m) (1 + n) / 2 \over N + (1 + m) (1 + n) / 2} + \sum_{i = 0}^{m} p_i E_{\p } \Big[ \log {X_i + (1 + n)  / 2 \over X_i + Z_i + (1 + n)  / 2} \Big] \non \\
&= \log {N + 1 + (1 + m) (1 + n) / 2 \over N + (1 + m) (1 + n) / 2} + \sum_{i = 0}^{m} {p_i}^2 E_{\p } \Big[ \log {X_i + (1 + n)  / 2 \over X_i + 1 + (1 + n)  / 2} \Big] \text{,} \non 
\end{align}
where the second equality follows since $Z_i \sim {\rm{Ber}} ( p_i )$ for each $i = 0, 1, \dots m$. 
\hfill$\Box$

\bigskip

\noindent
{\bf Proof of Lemma \ref{lem:convexity}.} \ \ The first half of this proof is similar to the proof of Lemma 4.2 of Hamura (2021a). 
Let $h \colon [0, \infty ) \to \mathbb{R}$ and let $h(x) = 0$ for $x \in (- \infty , 0)$. 
Let $X = X(p)$ for notational simplicity. 
Then 
\begin{align}
{\pd \over \pd p} E[ h(X) ] &= E \Big[ \Big( {X \over p} - {N - X \over 1 - p} \Big) h(X) \Big] = {E[ (X - N p) h(X)) ] \over p (1 - p)} = {E[ X \{ h(X) - h(X - 1) \} ] \over p} \text{,} \non 
\end{align}
where the last equality follows from the lemma of Johnson (1987). 
Therefore, 
\begin{align}
{\pd ^2 \over ( \pd p)^2} \{ p^2 E[ h(X) ] \} &= {\pd \over \pd p} [2 p E[ h(X) ] + p E[ X \{ h(X) - h(X - 1) \} ]] \non \\
&= 2 E[ h(X) ] + 2 E[ X \{ h(X) - h(X - 1) \} ] + E[ X \{ h(X) - h(X - 1) \} ] \non \\
&\quad + E[ X [X \{ h(X) - h(X - 1) \} - (X - 1) \{ h(X - 1) - h(X - 2) \} ] ] \non \\
&= E[ ( X^2 + 3 X + 2) h(X) - (2 X^2 + 2 X) h(X - 1) + X (X - 1) h(X - 2) ] \non \\
&= E[ \tilde{h} (X) - 2 \tilde{h} (X - 1) + \tilde{h} (X - 2) ] \text{,} \non 
\end{align}
where $\tilde{h} (x) = (x + 2) (x + 1) h(x)$ for $x \in \mathbb{R}$. 

Now, let $\nt = (1 + n) / 2$ and suppose that $h(x) = - \log \{ 1 - 1 / (x + 1 + \nt ) \} $ for all $x \in [0, \infty )$. 
Fix $x \in [0, \infty )$. 
Then 
\begin{align}
\tilde{h} (x) &= ( x^2 + 3 x + 2) \log {x + 1 + \nt \over x + \nt } \text{,} \non \\
{\tilde{h}}' (x) &= (2 x + 3) \log {x + 1 + \nt \over x + \nt } - {(x + 2) (x + 1) \over (x + 1 + \nt ) ( x + \nt )} \text{,} \quad \text{and} \non \\
{\tilde{h}}'' (x) &= 2 \log {x + 1 + \nt \over x + \nt } - {2 x + 3 \over (x + 1 + \nt ) (x + \nt )} \times 2 + {(x + 2) (x + 1) \over (x + 1 + \nt )^2 (x + \nt )} + {(x + 2) (x + 1) \over (x + 1 + \nt ) (x + \nt )^2} \text{.} \non 
\end{align}
Therefore, 
\begin{align}
( x + 1 + \nt ) {\tilde{h}}'' (x) &= 2 \sum_{k = 1}^{\infty } {1 \over k} {1 \over (x + 1 + \nt )^{k - 1}} - {2 x + 3 \over x + \nt } \times 2 + {(x + 2) (x + 1) \over (x + 1 + \nt ) (x + \nt )} + {(x + 2) (x + 1) \over (x + \nt )^2} \non \\
&= 2 + 2 \sum_{k = 2}^{\infty } {1 \over k} {1 \over (x + 1 + \nt )^{k - 1}} - 4 \Big( 1 - {\nt - 3 / 2 \over x + \nt } \Big) \non \\
&\quad + \Big( 1 - {\nt - 1 \over x + 1 + \nt } + 1 - {\nt - 2 \over x + \nt } \Big) \Big( 1 - {\nt - 1 \over x + \nt } \Big) \non \\
&= 2 \sum_{k = 2}^{\infty } {1 \over k} {1 \over (x + 1 + \nt )^{k - 1}} + 4 {\nt - 3 / 2 \over x + \nt } \non \\
&\quad - \Big( {\nt - 1 \over x + 1 + \nt } + {\nt - 2 \over x + \nt } \Big) - 2 {\nt - 1 \over x + \nt } + \Big( {\nt - 1 \over x + 1 + \nt } + {\nt - 2 \over x + \nt } \Big) {\nt - 1 \over x + \nt } \non 
\end{align}
and we have 
\begin{align}
( x + 1 + \nt ) {\tilde{h}}'' (x) &= 2 \sum_{k = 2}^{\infty } {1 \over k} {1 \over (x + 1 + \nt )^{k - 1}} + 4 {\nt - 3 / 2 \over x + \nt } \non \\
&\quad - \Big( {\nt - 1 \over x + \nt } - {\nt - 1 \over (x + \nt ) (x + 1 + \nt )} + {\nt - 2 \over x + \nt } \Big) - 2 {\nt - 1 \over x + \nt } + \Big( {\nt - 1 \over x + 1 + \nt } + {\nt - 2 \over x + \nt } \Big) {\nt - 1 \over x + \nt } \non \\
&= 2 \sum_{k = 2}^{\infty } {1 \over k} {1 \over (x + 1 + \nt )^{k - 1}} + {- 1 \over x + \nt } + {\nt - 1 \over (x + \nt ) (x + 1 + \nt )} + \Big( {\nt - 1 \over x + 1 + \nt } + {\nt - 2 \over x + \nt } \Big) {\nt - 1 \over x + \nt } \text{.} \non 
\end{align}
Since $\nt \ge 2$ by assumption, it follows that 
\begin{align}
( x + 1 + \nt ) {\tilde{h}}'' (x) &\ge 2 {1 \over 2} {1 \over x + 1 + \nt } + {- 1 \over x + \nt } + {\nt - 1 \over (x + \nt ) (x + 1 + \nt )} \non \\
&= - {1 \over (x + \nt ) (x + 1 + \nt )} + {\nt - 1 \over (x + \nt ) (x + 1 + \nt )} \ge 0 \text{.} \non 
\end{align}
Thus, $\tilde{h} (X) - 2 \tilde{h} (X - 1) + \tilde{h} (X - 2) \ge 0$ when $X \ge 2$. 
Also, 
\begin{align}
{\tilde{h} (1) - 2 \tilde{h} (0) + \tilde{h} (- 1) \over 6} &= \log \Big( 1 + {1 \over 1 + \nt } \Big) - {2 \over 3} \log \Big( 1 + {1 \over \nt } \Big) \ge \log \Big( 1 + {1 \over 1 + \nt } \Big) - \log \Big( 1 + {2 / 3 \over \nt } \Big) \ge 0 \text{.} \non 
\end{align}
It is trivial that $\tilde{h} (0) - 2 \tilde{h} (- 1) + \tilde{h} (- 2) \ge 0$. 
Hence, $\{ \pd ^2 / ( \pd p)^2 \} \{ p^2 E[ h(X) ] \} = E[ \tilde{h} (X) - 2 \tilde{h} (X - 1) + \tilde{h} (X - 2) ] \ge 0$ and this completes the proof. 
\hfill$\Box$

\bigskip

\noindent
{\bf Proof of Lemma \ref{lem:max_RD}.} \ \ It follows from Lemmas \ref{lem:RD} and \ref{lem:convexity} that 
\begin{align}
\De _{\p } (N, 1) &= \log \Big\{ 1 + {1 \over N + (1 + m) (1 + n) / 2} \Big\} - (1 + m) {1 \over 1 + m} \sum_{i = 0}^{m} H( p_i ) \non \\
&\le \log \Big\{ 1 + {1 \over N + (1 + m) (1 + n) / 2} \Big\} - (1 + m) H \Big( {1 \over 1 + m} \sum_{i = 0}^{m} p_i \Big) \non \\
&= \log \Big\{ 1 + {1 \over N + (1 + m) (1 + n) / 2} \Big\} - \sum_{i = 0}^{m} H \Big( {1 \over 1 + m} \Big) = \De _{\p ^{*}} (N, 1) \non 
\end{align}
by Jensen's inequality. 
\hfill$\Box$

\bigskip

\noindent
{\bf Proof of Theorem \ref{thm:N_Y}.} \ \ By Lemmas \ref{lem:decomposition} and \ref{lem:max_RD}, 
\begin{align}
\sup_{\p \in \Th } \De _{\p } (N, N' ) &= \sup_{\p \in \Th } \sum_{\nu = 1}^{N'} \De _{\p } (N + \nu - 1, 1) \non \\
&\le \sum_{\nu = 1}^{N'} \sup_{\p \in \Th } \De _{\p } (N + \nu - 1, 1) \non \\
&= \sum_{\nu = 1}^{N'} \De _{\p ^{*}} (N + \nu - 1, 1) = \De _{\p ^{*}} (N, N' ) \text{.} \non 
\end{align}
By Lemma \ref{lem:decomposition}, 
\begin{align}
\De _{\p ^{*}} (N, N' ) &= E_{\p ^{*}} \Big[ \log {N + N' + (1 + m) (1 + n) / 2 \over N + (1 + m) (1 + n) / 2} + \sum_{i = 0}^{m} {1 \over 1 + m} \log {X_i + (1 + n) / 2 \over X_i + Y_i + (1 + n) / 2} \Big] \non \\
&= - E_{\p ^{*}} \Big[ \sum_{i = 0}^{m} {1 \over 1 + m} \log \Big[ {Y_i / N' + \{ X_i + (1 + n) / 2 \} / N' \over 1 + \{ N + (1 + m) (1 + n) / 2 \} / N'} / {X_i + (1 + n) / 2 \over N + (1 + m) (1 + n) / 2} \Big] \Big] \non \\
&\to - E_{\p ^{*}} \Big[ \sum_{i = 0}^{m} {1 \over 1 + m} \log \Big\{ {1 \over 1 + m} / {X_i + (1 + n) / 2 \over N + (1 + m) (1 + n) / 2} \Big\} \Big] \non 
\end{align}
as $N' \to \infty $ by the law of large numbers. 
The right-hand side is negative since 
\begin{align}
P_{\p ^{*}} \Big( \Big( {X_i + (1 + n) / 2 \over N + (1 + m) (1 + n) / 2} \Big) _{i = 0}^{m} \neq \Big( {1 \over 1 + m} \Big) _{i = 0}^{m} \Big) > 0 \text{.} \non 
\end{align}
Thus, $\sup_{\p \in \Th } \De _{\p } (N, N' ) \le \De _{\p ^{*}} (N, N' ) < 0$ for sufficiently large $N'$. 
\hfill$\Box$

\bigskip

\noindent
{\bf Proof of Theorem \ref{thm:N_X}.} \ \ By Lemma \ref{lem:decomposition}, we can assume that $N' = 1$ without loss of generality. 
Then, by Lemmas \ref{lem:RD} and \ref{lem:max_RD}, it is sufficient to show that 
\begin{align}
\De _{\p ^{*}} (N, 1) &= \log \Big\{ 1 + {1 \over N + (1 + m) \nt } \Big\} - (1 + m) H \Big( {1 \over 1 + m} \Big) \non 
\end{align}
is negative when $N$ is sufficiently large, where $\nt = (1 + n) / 2$. 

Let $p^{*} = 1 / (1 + m)$ and $X^{*} \sim {\rm{Bin}} (N, p^{*} )$. 
Then 
\begin{align}
(1 + m) H \Big( {1 \over 1 + m} \Big) &= p^{*} E \Big[ - \log \Big( 1 - {1 \over X + 1 + \nt } \Big) \Big] = p^{*} \sum_{k = 1}^{\infty } {1 \over k} E \Big[ {1 \over (X + 1 + \nt )^k} \Big] \text{.} \non 
\end{align}
Note that, as in (3) of Cressie et al. (1981), 
\begin{align}
{1 \over k} E \Big[ {1 \over (X + 1 + \nt )^k} \Big] &= {1 \over k} {1 \over \Ga (k)} \int_{0}^{\infty } t^{k - 1} E[ e^{- (X + 1 + \nt ) t} ] dt = {1 \over k !} \int_{0}^{\infty } t^{k - 1} e^{- (1 + \nt ) t} E[ e^{- X t} ] dt \non 
\end{align}
for all $k = 1, 2, \dotsc $. 
Then 
\begin{align}
(1 + m) H \Big( {1 \over 1 + m} \Big) &= p^{*} \int_{0}^{\infty } \sum_{k = 1}^{\infty } {t^{k - 1} \over k !} e^{- (1 + \nt ) t} E[ e^{- X t} ] dt = p^{*} \int_{0}^{\infty } {e^t - 1 \over t} e^{- (1 + \nt ) t} E[ e^{- X t} ] dt \text{.} \non 
\end{align}
Since $E[ e^{- X t} ] = (1 - p^{*} + p^{*} e^{- t} )^N$ for all $t \in (0, \infty )$ (see, for example, Section 3.3 of Johnson, Kemp and Kotz (2005)), we have 
\begin{align}
(1 + m) H \Big( {1 \over 1 + m} \Big) &= p^{*} \int_{0}^{\infty } {e^t - 1 \over t} e^{- (1 + \nt ) t} (1 - p^{*} + p^{*} e^{- t} )^N dt \text{.} \non 
\end{align}
Meanwhile, 
\begin{align}
\log \Big\{ 1 + {1 \over N + (1 + m) \nt } \Big\} &= - \log \Big( 1 - {1 \over N + 1 + \nt / p^{*}} \Big) = \sum_{k = 1}^{\infty } {1 \over k} {1 \over (N + 1 + \nt / p^{*} )^k} \non \\
&= \sum_{k = 1}^{\infty } {1 \over k !} \int_{0}^{\infty } t^{k - 1} e^{- (N + 1 + \nt / p^{*} ) t} dt = \int_{0}^{\infty } {e^t - 1 \over t} e^{- (N + 1) t} e^{- \nt t / p^{*}} dt \text{.} \non 
\end{align}
Therefore, 
\begin{align}
\De _{\p ^{*}} (N, 1) &= \int_{0}^{\infty } {e^t - 1 \over t} e^{- (N + 1) t} e^{- \nt t / p^{*}} dt - p^{*} \int_{0}^{\infty } {e^t - 1 \over t} e^{- (1 + \nt ) t} (1 - p^{*} + p^{*} e^{- t} )^N dt \text{.} \non 
\end{align}

Note that 
\begin{align}
&0 \le {e^t - 1 \over t} = \sum_{k = 1}^{\infty } {t^{k - 1} \over k !} \le \sum_{k = 1}^{\infty } {t^{k - 1} \over (k - 1) !} = e^t \text{,} \non \\
&0 \le \Big( {e^t - 1 \over t} \Big) ' = \sum_{k = 2}^{\infty } {(k - 1) t^{k - 2} \over k !} \le \sum_{k = 2}^{\infty } {t^{k - 2} \over (k - 2) !} = e^t \text{,} \quad \text{and} \non \\
&0 \le \Big( {e^t - 1 \over t} \Big) '' = \sum_{k = 3}^{\infty } {(k - 1) (k - 2) t^{k - 3} \over k !} \le \sum_{k = 3}^{\infty } {t^{k - 3} \over (k - 3) !} = e^t \non 
\end{align}
and that 
\begin{align}
&{e^t - 1 \over t} \le {e^t \over t} \text{,} \non \\
&\Big( {e^t - 1 \over t} \Big) ' = {e^t \over t} - {e^t - 1\over t^2} \le {e^t \over t} \text{,} \quad \text{and} \non \\
&\Big( {e^t - 1 \over t} \Big) '' = {e^t \over t} - {e^t \over t^2} - {e^t \over t^2} + 2 {e^t - 1 \over t^3} \le {e^t \over t} + 2 {e^t \over t^3} \text{.} \non 
\end{align}
Also, note that $\nt \ge 2$ by assumption. 
Then, by integration by parts, 
\begin{align}
(N + 1) \De _{\p ^{*}} (N, 1) &= \Big[ - {e^t - 1 \over t} e^{- (N + 1) t} e^{- \nt t / p^{*}} \Big] _{0}^{\infty } - \Big[ - {e^t - 1 \over t} e^{- \nt t} (1 - p^{*} + p^{*} e^{- t} )^{N + 1} \Big] _{0}^{\infty } \non \\
&\quad + \int_{0}^{\infty } \Big\{ \Big( {e^t - 1 \over t} \Big) ' - {e^t - 1 \over t} {\nt \over p^{*}} \Big\} e^{- (N + 1) t} e^{- \nt t / p^{*}} dt \non \\
&\quad - \int_{0}^{\infty } \Big\{ \Big( {e^t - 1 \over t} \Big) ' - {e^t - 1 \over t} \nt \Big\} e^{- \nt t} (1 - p^{*} + p^{*} e^{- t} )^{N + 1} dt \non \\
&= \int_{0}^{\infty } g_1 (t) e^{- (N + 1) t} e^{- \nt t / p^{*}} dt - \int_{0}^{\infty } g_2 (t) e^{- \nt t} (1 - p^{*} + p^{*} e^{- t} )^{N + 1} dt \non 
\end{align}
where 
\begin{align}
g_1 (t) &= \Big( {e^t - 1 \over t} \Big) ' - {e^t - 1 \over t} {\nt \over p^{*}} \quad \text{and} \quad g_2 (t) = \Big( {e^t - 1 \over t} \Big) ' - {e^t - 1 \over t} \nt \non 
\end{align}
for $t \in (0, \infty )$. 
Furthermore, by integration by parts, 
\begin{align}
(N + 1)^2 \De _{\p ^{*}} (N, 1) &= \Big[ - g_1 (t) e^{- (N + 1) t} e^{- \nt t / p^{*}} \Big] _{0}^{\infty } \non \\
&\quad - {N + 1 \over N + 2} \Big[ - {1 \over p^{*}} g_2 (t) e^{- ( \nt - 1) t} (1 - p^{*} + p^{*} e^{- t} )^{N + 2} \Big] _{0}^{\infty } \non \\
&\quad + \int_{0}^{\infty } \Big\{ {g_1}' (t) - g_1 (t) {\nt \over p^{*}} \Big\} e^{- (N + 1) t} e^{- \nt t / p^{*}} dt \non \\
&\quad - {N + 1 \over N + 2} \int_{0}^{\infty } {1 \over p^{*}} \{ {g_2}' (t) - g_2 (t) ( \nt - 1) \} e^{- ( \nt - 1) t} (1 - p^{*} + p^{*} e^{- t} )^{N + 2} dt \non \\
&= g_1 (0) - {N + 1 \over N + 2} {1 \over p^{*}} g_2 (0) + \int_{0}^{\infty } \Big\{ {g_1}' (t) - g_1 (t) {\nt \over p^{*}} \Big\} e^{- (N + 1) t} e^{- \nt t / p^{*}} dt \non \\
&\quad - {N + 1 \over N + 2} \int_{0}^{\infty } {1 \over p^{*}} \{ {g_2}' (t) - g_2 (t) ( \nt - 1) \} e^{- ( \nt - 1) t} (1 - p^{*} + p^{*} e^{- t} )^{N + 2} dt \text{,} \non 
\end{align}
since 
\begin{align}
| \{ {g_2}' (t) - g_2 (t) ( \nt - 1) \} e^{- ( \nt - 1) t} | &\le \{ e^t + e^t \nt ( \nt - 1) \} e^{- ( \nt - 1) t} = \{ 1 + \nt ( \nt - 1) \} e^{- ( \nt - 2) t} \non 
\end{align}
if $\nt > 2$ and since 
\begin{align}
| \{ {g_2}' (t) - g_2 (t) ( \nt - 1) \} e^{- ( \nt - 1) t} | &\begin{cases} \displaystyle %
\le | {g_2}' (t) - g_2 (t)| \le 6 e^t\text{,}  & \text{for all $t \in (0, 1]$} \text{,} \\ \displaystyle = \Big| 2 {e^t \over t^3} + {e^t \over t^2} - {2 \over t^3} - {3 \over t^2} - {2 \over t} \Big| e^{- t} \le {3 \over t^2} + {7 \over t} e^{- t} \text{,} & \text{for all $t \in (1, \infty )$} \text{,} \end{cases} \non 
\end{align}
if $\nt = 2$. 
Thus, by the dominated convergence theorem, 
\begin{align}
(N + 1)^2 \De _{\p ^{*}} (N, 1) &\to g_1 (0) - {1 \over p^{*}} g_2 (0) = {1 \over 2} \Big( 1 - {1 \over p^{*}} \Big) < 0 \non 
\end{align}
as $N \to \infty $, and this completes the proof. 
\hfill$\Box$


\begin{thebibliography}{00}



\bibitem{am1986}
Alam, K. and Mitra, A. (1986). 
An empirical bayes estimate of multinomial probabilities. 
{\it Communications in Statistics - Theory and Methods}, {\bf 15}, 3103--3127. 


\bibitem{a1981}
Albert, J. (1981). 
Pseudo-bayes estimation of multinomial proportions. 
{\it Communications in Statistics - Theory and Methods}, {\bf 10}, 1587--1611. 



\bibitem{a1985}
Albert, J. (1985). 
Bayesian estimation methods for incomplete two-way contingency tables using prior beliefs of association. 
{\it Bayesian Statistics}, {\bf 2}. 



\bibitem{a1987}
Albert, J.H. (1987). 
Empirical bayes estimation in contingency tables. 
{\it Communications in Statistics - Theory and Methods}, {\bf 16}, 2459--2485. 


\bibitem{ag1983}
Albert, J.H. and Gupta, A.K. (1983). 
Bayesian estimation methods for $2 \times 2$ contingency tables using mixtures of Dirichlet distributions. 
{\it Journal of the American Statistical Association}, {\bf 78}, 708--717. 


\bibitem{cdfp1981}
Cressie, N., Davis, A.S., Folks, J.L. and Policello II, G.E. (1981). 
Moment-generating function and negative integer foments. 
{\it The American Statistician}, {\bf 35}, 148--150. 


\bibitem{fh1973}
Fienberg, S.E. and Holland, P.W. (1973). 
Simultaneous estimation of multinomial cell probabilities. 
{\it Journal of the American Statistical Association}, {\bf 68}, 683--691. 


\bibitem{gg1994}
Gibbons, P.C. and Greenberg, E. (1994). 
Bayesian reconstruction of contingency tables with partially categorized data. 
{\it Communications in Statistics - Theory and Methods}, {\bf 23}, 3349--3359. 



\bibitem{h2021a}
Hamura, Y. (2021a). 
Bayesian point estimation and predictive density estimation for the binomial distribution with a restricted probability parameter. 
{\it Communications in
Statistics - Theory and Methods}. 

\bibitem{h2021b}
Hamura, Y. (2021b). 
Bayesian shrinkage estimation for stratified count data. 
{\it arXiv preprint  arXiv:2112.13245}. 





\bibitem{j1987}
Johnson, R.W. (1987). 
Simultaneous estimation of binomial $N$'s. 
{\it Sankhya Series A}, {\bf 49}, 264--267. 



\bibitem{jkk2005}
Johnson, N.L., Kemp, A.W. and Kotz, S. (2005). 
Univariate discrete distributions. 
Wiley,New Jersey, 2005. 



\bibitem{k2012}
Komaki, F. (2012). 
Asymptotically minimax Bayesian predictive densities for multinomial models. 
{\it Electronic Journal of Statistics}, {\bf 6}, 934--957. 




\bibitem{l1977}
Leonard, T. (1977). 
A Bayesian approach to some multinomial estimation and pretesting problems. 
{\it Journal of the American Statistical Association}, {\bf 72}, 869--874. 



\bibitem{nn1989}
Nayak, T.K. and Naik, D.N. (1989). 
Estimating multinomial cell probabilities under quadratic loss. 
{\it Journal of the Royal Statistical Society. Series D (The Statistician)}, {\bf 38}, 3--10. 








\end{thebibliography}
\end{document}